\documentclass[12pt,dvips]{amsart}
\usepackage{euler, amsfonts, amssymb, latexsym, epsfig,epic}

\newtheorem{Theorem}{Theorem}

\newtheorem{Lemma}{Lemma}

\newtheorem{Corollary}{Corollary}
\newtheorem*{Corollary*}{Corollary}
 
\newtheorem*{Theorem*}{Theorem}

\theoremstyle{remark}
\newtheorem{Example}{Example}

\newcommand\complexes{{\mathbb C}}

\newcommand\naturals{{\mathbb N}}

\theoremstyle{plain}

\newcommand\dfn{\bf} 

\newcommand\Fln{{{\rm Flags}(\complexes^n)}}
\newcommand\groth{{\mathfrak G}}

\font\co=lcircle10

\def\jr{\smash{\raise2pt\hbox{\co \rlap{\rlap{\char'005} \char'007}}
               \raise6pt\hbox{\rlap{\vrule height5pt}}
               \raise2pt\hbox{\rlap{\hskip4pt \vrule height0.4pt depth0pt
                width5.7pt}}
               \raise2pt\hbox{\rlap{\hskip-9.5pt \vrule height.4pt depth0pt
                width6.2pt}}
               \lower6pt\hbox{\rlap{\vrule height4.5pt}}}}
\def\je{\smash{\raise2pt\hbox{\co \rlap{\rlap{\char'005}
                \phantom{\char'007}}}\raise6pt\hbox{\rlap{\vrule height5pt}}
               \raise2pt\hbox{\rlap{\hskip-9.5pt \vrule height.4pt depth0pt
                width6.2pt}}}}
\def\er{\smash{\raise2pt\hbox{\co \rlap{\rlap{\phantom{\char'005}} \char'007}}
               \raise2pt\hbox{\rlap{\hskip4pt \vrule height0.4pt depth0pt
                width5.7pt}}
               \lower6pt\hbox{\rlap{\vrule height4.5pt}}}}
\def\+{\smash{\lower6pt\hbox{\rlap{\vrule height17pt}}
                \raise2pt\hbox{\rlap{\hskip-9pt \vrule height.4pt depth0pt
                width18.7pt}}}}
\def\hor{\smash{\raise2pt\hbox{\rlap{\hskip-9.5pt \vrule height.4pt depth0pt
                width19.2pt}}}}
\def\ver{\smash{\lower6pt\hbox{\rlap{\vrule height17pt}}}}

\def\textcross{\ \smash{\lower4pt\hbox{\rlap{\hskip4.15pt\vrule height14pt}}
                \raise2.8pt\hbox{\rlap{\hskip-3pt \vrule height.4pt depth0pt
                width14.7pt}}}\hskip12.7pt}
\def\textelbow{\ \hskip.1pt\smash{\raise2.8pt%
                \hbox{\co \hskip 4.15pt\rlap{\rlap{\char'005} \char'007}
                \lower6.8pt\rlap{\vrule height3.5pt}
                \raise3.6pt\rlap{\vrule height3.5pt}}
                \raise2.8pt\hbox{%
                  \rlap{\hskip-7.15pt \vrule height.4pt depth0pt width3.5pt}%
                  \rlap{\hskip4.05pt \vrule height.4pt depth0pt width3.5pt}}}
                \hskip8.7pt}

\newcommand{\cellsize}{22}
\newlength{\cellsz} 
\newsavebox{\cell}
\sbox{\cell}{\begin{picture}(\cellsize,\cellsize)
\put(0,0){\line(1,0){\cellsize}}
\put(0,0){\line(0,1){\cellsize}}
\put(\cellsize,0){\line(0,1){\cellsize}}
\put(0,\cellsize){\line(1,0){\cellsize}}
\end{picture}}
\newcommand\cellify[1]{\def\thearg{#1}\def\nothing{}%
\ifx\thearg\nothing
\vrule width0pt height\cellsz depth0pt\else
\hbox to 0pt{\usebox{\cell} \hss}\fi%
\vbox to \cellsz{
\vss
\hbox to \cellsz{\hss$#1$\hss}
\vss}}
\newcommand\tableau[1]{\vtop{\let\\\cr
\baselineskip -16000pt \lineskiplimit 16000pt \lineskip 0pt
\ialign{&\cellify{##}\cr#1\crcr}}}
\begin{document}
\pagestyle{plain}

\title{A formula for $K$-theory truncation Schubert calculus}
\author{Allen Knutson}
\address{Department of Mathematics, University of California, Berkeley, CA 94720, USA}
\email{allenk@math.berkeley.edu}
\thanks{AK was supported by an NSF grant.}

\author{Alexander Yong}
\address{Department of Mathematics, University of California, Berkeley, CA 94720, USA}
\email{ayong@math.berkeley.edu}
\date{Draft of \today}

\maketitle

\begin{abstract}
Define a {\dfn truncation} $r_t(p)$ of a polynomial $p$ in 
$\{x_1,x_2,x_3,\ldots\}$ as the polynomial with all but the first $t$ variables
set to zero. In certain good cases, the truncation of a Schubert or
Grothendieck polynomial may again be a Schubert or Grothendieck polynomial.
We use this phenomenon to give subtraction-free formulae for certain Schubert 
structure constants in $K(\Fln)$, in particular 
generalizing those from [Kogan, '00] in which only cohomology was treated,
and from [Buch, '02] on the Grassmannian case. The terms in the answer
are computed using ``marching'' operations on permutation diagrams.
\end{abstract}


\section{Introduction}

Let $\Fln$ denote the variety of
complete flags in ${\mathbb C}^n$. To each permutation $\pi$
in the symmetric group $S_{n}$, 
there is an associated Schubert variety $X_{\pi} \subseteq \Fln$. 
The classes of the Schubert structure sheaves $[{\mathcal O}_{X_{\pi}}]$ form
an additive ${\mathbb Z}$-linear basis of the $K$-theory (Grothendieck) ring 
$K(\Fln)$ of algebraic vector bundles over $\Fln$. The
{\dfn Schubert structure constants} are the integers defined by
\begin{equation}
\label{eqn:main_def}
[{\mathcal O}_{X_{\sigma}}]\cdot [{\mathcal O}_{X_{\rho}}]
=\sum_{\pi\in S_{n}} C_{\sigma,\rho}^{\pi} [{\mathcal O}_{X_{\pi}}].
\end{equation}
It is known~\cite{brion:pos} that 
$(-1)^{\ell(\sigma)+\ell(\rho)-\ell(\pi)}C_{\sigma,\rho}^{\pi}\geq 0$, 
where $\ell(\alpha)$ is the minimum $\ell$ such that $\alpha$ is expressible as
a product of $\ell$ simple transpositions $s_i = t_{i\leftrightarrow i+1}$. 
In the {\em cohomology case}, i.e., when 
$\ell(\sigma)+\ell(\rho)=\ell(\pi)$, these are 
the structure constants for the analogous expansion of the product of 
Schubert classes $[X_{\sigma}]\cdot[X_{\rho}]$  
in the cohomology ring $H^{*}(\Fln)$,
and 
count the number of points in the
intersection of general triple translates of $X_{\sigma}, X_{\rho}$ and
$X_{w_0 \pi}$ (where $w_0$ denotes the longest permutation in $S_n$).
The expansion (\ref{eqn:main_def}) behaves well with respect to the 
inclusion $S_n \hookrightarrow S_{n+1}$. In particular, for 
any two permutations $\sigma,\rho$, and $n$ sufficiently large, 
(\ref{eqn:main_def}) stabilizes. Therefore, it will be unambiguous
(and convenient) to call $(\sigma,\rho,\pi)\in S_{\infty}^3$ a 
{\dfn Schubert problem}, where $S_{\infty}=\bigcup_{n\geq 1} S_n$.

It is a famous open problem to give a general subtraction-free
combinatorial formula applicable to any Schubert problem.  The
analogous problem for Grassmannians is solved in the cohomology case
by the Littlewood-Richardson rule, and more recently in $K$-theory
by A.~Buch~\cite{buch:KLR}. A solution for the flag variety would provide 
an important generalization of the Littlewood-Richardson rule. 
However, the known generalized Littlewood-Richardson rules
handle only limited cases of the Schubert problems,
with \cite{kogan:schur_schub,buch:KLR} representing the state of the art.

Our main result is a subtraction-free combinatorial 
formula for the family we call {\em truncation Schubert problems}
(defined below). 
This formula specializes to compute the 
$K$-theory generalizations of the numbers considered by 
M.~Kogan~\cite{kogan:schur_schub} and the $K$-theory \linebreak 
Littlewood-Richardson coefficients of~\cite{buch:KLR}.
Actually, our main result gives formulas
for many other combinatorial
numbers studied in connection to the Schubert 
calculus~\cite{buch:KLR, lascoux, bergeron.sottile:duke}, 
formulas for Schubert and Grothendieck
polynomials~\cite{lascoux.schutzenberger:hopf, BJS, fomin.kirillov:1,
fomin.kirillov:2}, 
degeneracy loci~\cite{buch.fulton:chern, buch:duke, fulton:univ, BKTY:1,
BKTY:2} and quantum Schubert polynomials~\cite{FGP, CF:99, BKTY:1},
see, e.g.,~\cite{buch.sottile.yong:quiver} and the references therein. We 
find it interesting that our
formula also applies to new cases of Schubert problems where {\em neither}
class is a pullback from a Grassmannian.

The fact that these numbers all arise from the 
single framework of problems isolated here suggests  
that their common 
combinatorial and geometric features ought to be better understood.
 
On the combinatorial side, in Section~3, we present our formula
in terms of simple ``marching'' moves of the diagram of a permutation. 
We would like to understand how, e.g., various combinatorial aspects
of the classical 
Littlewood-Richardson coefficients \linebreak might extend to this 
family of numbers. 
It would be interesting to understand the relations
between the formula given here and formulas for the aforementioned special
cases, and other related formulas, e.g.,~\cite{fulton.lascoux,lenart,
lenart.sottile:pieri,pittie.ram,sottile:pieri}. 

One feature of our proof is that it is both short and 
completely 
combinatorial. It is based on
``truncation'' techniques concerning Grothendieck 
polynomials~\cite{lascoux.schutzenberger:hopf} and in particular, the
``transition'' formula of A.~Lascoux~\cite{lascoux}. 
These methods (at least in cohomology) can be considered classical in the 
subject. Indeed, in
previous work~\cite{lascoux} (see also~\cite{lascoux.schutzenberger:LR}), 
similar techniques were applied to give new formulas for the 
$K$-theory Littlewood-Richardson coefficients (after~\cite{buch:KLR}).
However, it is perhaps surprising that such methods are in fact applicable
to more general Schubert calculus problems, and in particular, Kogan's 
Schubert problems. Our principal novelty of three simultaneous 
observations is reflected respectively 
in the three equalities found in (\ref{eqn:break}) from
Section~3. 

	Thus, since Kogan's Schubert problems form a special case of the 
truncation Schubert problems, our formula covers 
new cases of the Schubert problem in $K$-theory (and moreover, our proof 
makes transparent the role of Kogan's conditions). However, we emphasize 
that our formula handles new cases beyond that 
in~\cite{kogan:schur_schub}, even in cohomology.

A further goal of this paper is to present the diagram marching moves.
One reason to use such (recursive)
combinatorics is that the moves have a natural 
{\em geometric} interpretation. In a sequel~\cite{KYI} to this paper, 
we interpret
the moves in terms of Gr\"{o}bner degeneration 
of matrix Schubert varieties~\cite{fulton:duke} via {\em diagonal} term orders
(in an important contrast to the {\em anti}-diagonal term orders used
in \cite{KM:annals}). For example, in the cohomology case, 
our formula can be interpreted as counting certain
components of a partially degenerated matrix Schubert variety. 
It would be interesting to understand what 
relations exist between the formula presented here and the 
geometric Littlewood-Richardson rule of R.~Vakil~\cite{vakil:LR},
which is also based on degeneration.

Finally, one other advantage of the approach presented here is the 
possible extensions to other Schubert
calculus settings, e.g., the
cohomology/$K$-theory ring of flag varieties corresponding to the other 
classical Lie types $BCD$ (work in progress with F.~Sottile). 

We thank Frank Sottile and Alexander Woo for helpful discussions.

\section{Diagram moves and the main result}
Let $G(\pi)$ denote the permutation matrix 
associated to $\pi\in S_n$, and 
call the nonzero entries of $G(\pi)$ its {\dfn dots}. The
{\dfn diagram} of a permutation $\pi$ 
is the following subset of $[n]\times [n]$: 
\[D(\pi) := \{(p,q), \ 1\leq p,q\leq n, \ \pi(p)>q, \ \pi^{-1}(q)>p\}.\]
Equivalently, $(p,q) \in D(\pi)$ if $(p,\pi^{-1}(q))$ is an inversion
of $\pi$; thus $\#D(\pi) = \ell(\pi)$.
Graphically, $D(\pi)$ is obtained from $G(\pi)$ by drawing 
a ``hook'' consisting of lines going east and south from each dot.
The diagram appears as a collection of  ``connected components'' 
of squares not in the hook of any dot (see Example~1 below).

Call the southernmost, then eastmost, box $(l,m)\in D(\pi)$ the 
{\dfn maximal corner}. Note that the maximal corner of $\pi$
is in row $l$ if and only if the last {\dfn descent} of $\pi$ is in row $l$, 
i.e., the largest index $l$ such that $\pi(l)>\pi(l+1)$. 
Call any dot that is maximally southeast with respect to the condition
that it is northwest of $(l,m)$ a {\dfn pivot}. There are no pivots 
if and only if the maximal corner is in the connected
component of $D(\pi)$ attached to the top left corner of $[n]\times [n]$. 
See Example~1 below.

In the following definitions, 
it is convenient to describe the cohomology
versions first before explaining their $K$-theory analogues.
In the next section, we will connect what follows to the
Grothendieck transition formula of~\cite{lascoux}.

First, we describe the {\dfn marching} operation on $D(\pi)$.
Suppose the maximal corner is at $(l,m)$. If the input permutation has
no pivots, declare the output of the marching to be null
``$\emptyset$'', and write $\pi\to\emptyset$. Otherwise,
consider a pivot $(i,j)\in G(\pi)$. 
Remove the hook emanating from $(i,j)$, 
and move {\em strictly} to the northwest
every diagram box in the rectangle with the corners $(i,j),(l,m)$,
into the only spaces available (i.e., by 
``hopping'' over any hooks in the way). Do this by starting 
with the unique northwest box in the rectangle and continue left to right
along the rows, and from top to bottom. 
It is easy to check that the resulting collection of boxes is
necessarily the diagram of a permutation $\rho$. Let 
$\pi{\stackrel{i}{\longrightarrow}}\rho$ denote that $\rho$ 
is obtained from marching on $D(\pi)$ {\dfn towards the pivot} in row $i$.

\begin{Example}
Let $\pi=4317625$. We have
$D(\pi)=\{(1,1),(1,2),(1,3),(2,1),(2,2),(4,2),$ $(4,5),$ 
$(4,6),$ $(5,2),(5,5)\}.$
The maximal corner is $(5,5)$ and its pivots are the dots at
$(1,4),(2,3)$ and $(3,1)$.  
The boxes in the $(3,1),(5,5)$ rectangle of $\pi$ are marked with Xs.
Marching towards the pivot $(2,3)$ we get
$\rho=4517326$:
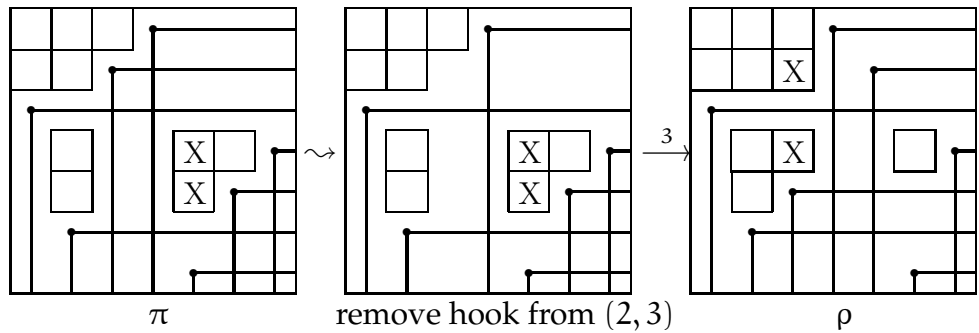
\setlength{\unitlength}{.27mm}
\[
\begin{picture}(465,155)
\put(68,0){$\pi$}
\put(406,0){$\rho$}
\put(160,0){remove hook from $(2,3)$}
\put(0,15){\framebox(140,140)}
\put(70,145){\circle*{4}}
\thicklines
\put(70,145){\line(1,0){70}}
\put(70,145){\line(0,-1){130}}
\put(50,125){\circle*{4}}
\put(50,125){\line(1,0){90}}
\put(50,125){\line(0,-1){110}}
\put(10,105){\circle*{4}}
\put(10,105){\line(1,0){130}}
\put(10,105){\line(0,-1){90}}
\put(130,85){\circle*{4}}
\put(130,85){\line(1,0){10}}
\put(130,85){\line(0,-1){70}}
\put(110,65){\circle*{4}}
\put(110,65){\line(1,0){30}}
\put(110,65){\line(0,-1){50}}
\put(30,45){\circle*{4}}
\put(30,45){\line(1,0){110}}
\put(30,45){\line(0,-1){30}}

\put(90,25){\circle*{4}}
\put(90,25){\line(1,0){50}}
\put(90,25){\line(0,-1){10}}

\thinlines

\put(0,135){\framebox(60,20)}
\put(20,135){\line(0,1){20}}
\put(40,135){\line(0,1){20}}
\put(0,115){\line(1,0){40}}
\put(40,115){\line(0,1){20}}
\put(20,115){\line(0,1){20}}
\put(20,55){\framebox(20,40)}
\put(85,59){{\rm X}}
\put(85,79){{\rm X}}
\put(20,75){\line(1,0){20}}
\put(80,55){\line(0,1){40}}
\put(80,55){\line(1,0){20}}
\put(100,55){\line(0,1){40}}
\put(80,95){\line(1,0){40}}
\put(80,75){\line(1,0){40}}
\put(120,75){\line(0,1){20}}


\put(250,59){{\rm X}}
\put(250,79){{\rm X}}
\put(165,15){\framebox(140,140)}
\put(235,145){\circle*{4}}
\thicklines
\put(235,145){\line(1,0){70}}
\put(235,145){\line(0,-1){130}}
\put(175,105){\circle*{4}}
\put(175,105){\line(1,0){130}}
\put(175,105){\line(0,-1){90}}
\put(295,85){\circle*{4}}
\put(295,85){\line(1,0){10}}
\put(295,85){\line(0,-1){70}}
\put(275,65){\circle*{4}}
\put(275,65){\line(1,0){30}}
\put(275,65){\line(0,-1){50}}
\put(195,45){\circle*{4}}
\put(195,45){\line(1,0){110}}
\put(195,45){\line(0,-1){30}}

\put(255,25){\circle*{4}}
\put(255,25){\line(1,0){50}}
\put(255,25){\line(0,-1){10}}

\thinlines

\put(165,135){\framebox(60,20)}
\put(185,135){\line(0,1){20}}
\put(205,135){\line(0,1){20}}
\put(165,115){\line(1,0){40}}
\put(205,115){\line(0,1){20}}
\put(185,115){\line(0,1){20}}
\put(185,55){\framebox(20,40)}
\put(185,75){\line(1,0){20}}
\put(245,55){\line(0,1){40}}
\put(245,55){\line(1,0){20}}
\put(265,55){\line(0,1){40}}
\put(245,95){\line(1,0){40}}
\put(245,75){\line(1,0){40}}
\put(285,75){\line(0,1){20}}

\put(145,80){$\leadsto$}
\put(310,80){$\stackrel{3}{\longrightarrow}$}

\put(380,119){{\rm X}}
\put(380,79){{\rm X}}

\put(335,15){\framebox(140,140)}
\put(405,145){\circle*{4}}
\thicklines
\put(405,145){\line(1,0){70}}
\put(405,145){\line(0,-1){130}}
\put(425,125){\circle*{4}}
\put(425,125){\line(1,0){50}}
\put(425,125){\line(0,-1){110}}
\put(345,105){\circle*{4}}
\put(345,105){\line(1,0){130}}
\put(345,105){\line(0,-1){90}}
\put(465,85){\circle*{4}}
\put(465,85){\line(1,0){10}}
\put(465,85){\line(0,-1){70}}

\put(385,65){\circle*{4}}
\put(385,65){\line(1,0){90}}
\put(385,65){\line(0,-1){50}}

\put(365,45){\circle*{4}}
\put(365,45){\line(1,0){110}}
\put(365,45){\line(0,-1){30}}

\put(445,25){\circle*{4}}
\put(445,25){\line(1,0){30}}
\put(445,25){\line(0,-1){10}}

\thinlines

\put(335,115){\framebox(60,40)}
\put(335,135){\line(1,0){60}}
\put(355,115){\line(0,1){40}}
\put(375,115){\line(0,1){40}}
\put(355,75){\framebox(40,20)}
\put(355,75){\line(0,-1){20}}
\put(355,55){\line(1,0){20}}
\put(375,55){\line(0,1){40}}
\put(435,75){\framebox(20,20)}
\end{picture}
\]
Note that two adjacent ${\rm X}$'s at $(4,5)$ and
$(5,5)$ can become separated (to $(2,3)$ and $(4,3)$ respectively)
after marching.

Marching instead towards the pivot $(3,1)$ we get $\rho=4357126$:
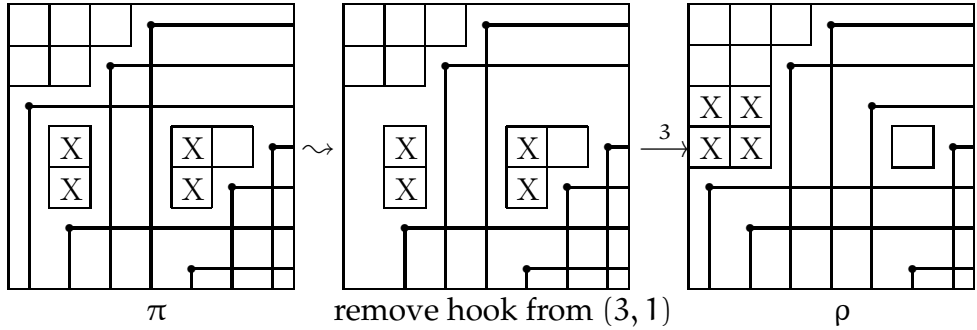
\setlength{\unitlength}{.27mm}
\[
\begin{picture}(465,155)
\put(68,0){$\pi$}
\put(406,0){$\rho$}
\put(160,0){remove hook from $(3,1)$}
\put(0,15){\framebox(140,140)}
\put(70,145){\circle*{4}}
\thicklines
\put(70,145){\line(1,0){70}}
\put(70,145){\line(0,-1){130}}
\put(50,125){\circle*{4}}
\put(50,125){\line(1,0){90}}
\put(50,125){\line(0,-1){110}}
\put(10,105){\circle*{4}}
\put(10,105){\line(1,0){130}}
\put(10,105){\line(0,-1){90}}
\put(130,85){\circle*{4}}
\put(130,85){\line(1,0){10}}
\put(130,85){\line(0,-1){70}}
\put(110,65){\circle*{4}}
\put(110,65){\line(1,0){30}}
\put(110,65){\line(0,-1){50}}
\put(30,45){\circle*{4}}
\put(30,45){\line(1,0){110}}
\put(30,45){\line(0,-1){30}}

\put(90,25){\circle*{4}}
\put(90,25){\line(1,0){50}}
\put(90,25){\line(0,-1){10}}

\thinlines

\put(0,135){\framebox(60,20)}
\put(20,135){\line(0,1){20}}
\put(40,135){\line(0,1){20}}
\put(0,115){\line(1,0){40}}
\put(40,115){\line(0,1){20}}
\put(20,115){\line(0,1){20}}
\put(20,55){\framebox(20,40)}
\put(25,59){\bf {\rm X}}
\put(25,79){\bf {\rm X}}
\put(85,59){{\rm X}}
\put(85,79){{\rm X}}
\put(20,75){\line(1,0){20}}
\put(80,55){\line(0,1){40}}
\put(80,55){\line(1,0){20}}
\put(100,55){\line(0,1){40}}
\put(80,95){\line(1,0){40}}
\put(80,75){\line(1,0){40}}
\put(120,75){\line(0,1){20}}


\put(190,59){\bf {\rm X}}
\put(190,79){\bf {\rm X}}
\put(250,59){{\rm X}}
\put(250,79){{\rm X}}
\put(165,15){\framebox(140,140)}
\put(235,145){\circle*{4}}
\thicklines
\put(235,145){\line(1,0){70}}
\put(235,145){\line(0,-1){130}}
\put(215,125){\circle*{4}}
\put(215,125){\line(1,0){90}}
\put(215,125){\line(0,-1){110}}
\put(295,85){\circle*{4}}
\put(295,85){\line(1,0){10}}
\put(295,85){\line(0,-1){70}}
\put(275,65){\circle*{4}}
\put(275,65){\line(1,0){30}}
\put(275,65){\line(0,-1){50}}
\put(195,45){\circle*{4}}
\put(195,45){\line(1,0){110}}
\put(195,45){\line(0,-1){30}}

\put(255,25){\circle*{4}}
\put(255,25){\line(1,0){50}}
\put(255,25){\line(0,-1){10}}

\thinlines

\put(165,135){\framebox(60,20)}
\put(185,135){\line(0,1){20}}
\put(205,135){\line(0,1){20}}
\put(165,115){\line(1,0){40}}
\put(205,115){\line(0,1){20}}
\put(185,115){\line(0,1){20}}
\put(185,55){\framebox(20,40)}
\put(185,75){\line(1,0){20}}
\put(245,55){\line(0,1){40}}
\put(245,55){\line(1,0){20}}
\put(265,55){\line(0,1){40}}
\put(245,95){\line(1,0){40}}
\put(245,75){\line(1,0){40}}
\put(285,75){\line(0,1){20}}

\put(145,80){$\leadsto$}
\put(310,80){$\stackrel{3}{\longrightarrow}$}

\put(340,99){\bf {\rm X}}
\put(340,79){\bf {\rm X}}
\put(360,99){{\rm X}}
\put(360,79){{\rm X}}

\put(335,15){\framebox(140,140)}
\put(405,145){\circle*{4}}
\thicklines
\put(405,145){\line(1,0){70}}
\put(405,145){\line(0,-1){130}}
\put(385,125){\circle*{4}}
\put(385,125){\line(1,0){90}}
\put(385,125){\line(0,-1){110}}
\put(425,105){\circle*{4}}
\put(425,105){\line(1,0){50}}
\put(425,105){\line(0,-1){90}}
\put(465,85){\circle*{4}}
\put(465,85){\line(1,0){10}}
\put(465,85){\line(0,-1){70}}

\put(345,65){\circle*{4}}
\put(345,65){\line(1,0){130}}
\put(345,65){\line(0,-1){50}}

\put(365,45){\circle*{4}}
\put(365,45){\line(1,0){110}}
\put(365,45){\line(0,-1){30}}

\put(445,25){\circle*{4}}
\put(445,25){\line(1,0){30}}
\put(445,25){\line(0,-1){10}}

\thinlines

\put(335,75){\framebox(40,80)}
\put(335,95){\line(1,0){40}}
\put(335,115){\line(1,0){40}}
\put(335,135){\line(1,0){60}}
\put(395,135){\line(0,1){20}}
\put(355,75){\line(0,1){80}}
\put(435,75){\framebox(20,20)}
\end{picture}
\]
This time, some nonadjacent boxes in $\pi$ become adjacent in $\rho$.
\end{Example}

More generally, suppose that $1\leq i_1<i_2<\ldots<i_k<l$ are a subset of 
the rows that contain pivots
of $(l,m)$. Consider $\pi\stackrel{i_1}{\longrightarrow}\rho_1$. 
Add a box to the diagram of $\rho_1$ at $(l,\rho_1(l))$. 
This is the diagram of a new permutation
$\rho_1'$, where the added box is the maximal corner for $\rho_1'$,
and the dot in row $i_2$ is a pivot for this box.
Now march $\rho_1'\stackrel{i_2}{\longrightarrow}\rho_2$, 
and add a maximal corner similarly in row $l$ to give $\rho_2'$. 
Repeat this process of marching and adding a box in row $l$ 
until we obtain $\rho=\rho_k$. 
We write $\pi\stackrel{i_1,i_2,\ldots,i_k}{\relbar\joinrel\relbar\joinrel\relbar\joinrel\relbar\joinrel\longrightarrow}\rho$ to 
denote this more general $K$-{\dfn marching} operation; a total of 
$k-1$ boxes are added.
\begin{Example}
Let $\pi$ be as in the above example and suppose we march to
the pivots $(1,4)$ and $(3,1)$ of the maximal corner 
$(5,5)$, in succession. 
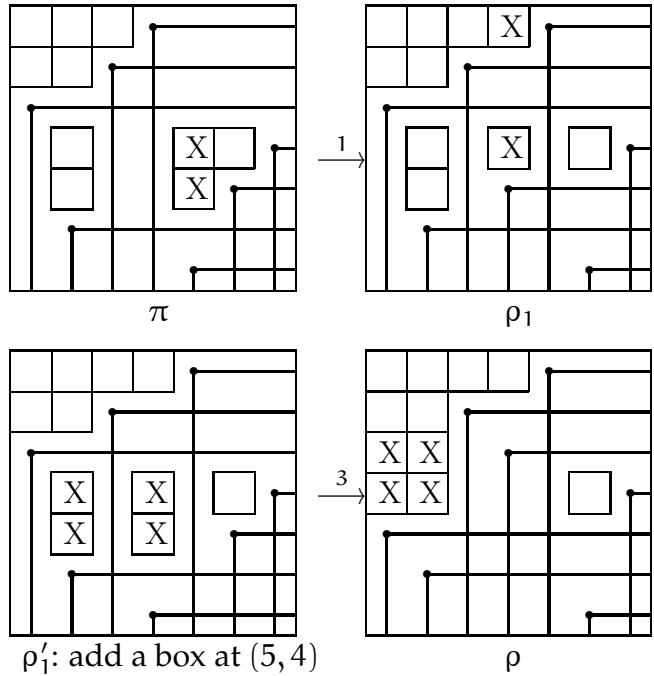
\setlength{\unitlength}{.27mm}
\[
\begin{picture}(330,325)
\put(68,170){$\pi$}
\put(243,170){$\rho_{1}$}
\put(5,0){$\rho_1'$: \mbox{add a box at $(5,4)$}}
\put(243,0){$\rho$}

\put(26,60){\rm X}
\put(26,80){\rm X}
\put(66,60){\rm X}
\put(66,80){\rm X}
\put(181,80){\rm X}
\put(201,80){\rm X}
\put(181,100){\rm X}
\put(201,100){\rm X}
\put(0,15){\framebox(140,140)}
\put(90,145){\circle*{4}}
\thicklines
\put(90,145){\line(1,0){50}}
\put(90,145){\line(0,-1){130}}
\put(50,125){\circle*{4}}
\put(50,125){\line(1,0){90}}
\put(50,125){\line(0,-1){110}}
\put(10,105){\circle*{4}}
\put(10,105){\line(1,0){130}}
\put(10,105){\line(0,-1){90}}
\put(130,85){\circle*{4}}
\put(130,85){\line(1,0){10}}
\put(130,85){\line(0,-1){70}}
\put(110,65){\circle*{4}}
\put(110,65){\line(1,0){30}}
\put(110,65){\line(0,-1){50}}
\put(30,45){\circle*{4}}
\put(30,45){\line(1,0){110}}
\put(30,45){\line(0,-1){30}}

\put(70,25){\circle*{4}}
\put(70,25){\line(1,0){70}}
\put(70,25){\line(0,-1){10}}

\thinlines

\put(0,135){\framebox(80,20)}
\put(20,135){\line(0,1){20}}
\put(40,135){\line(0,1){20}}
\put(60,135){\line(0,1){20}}
\put(0,115){\line(1,0){40}}
\put(40,115){\line(0,1){20}}
\put(20,115){\line(0,1){20}}
\put(20,55){\framebox(20,40)}
\put(20,75){\line(1,0){20}}
\put(60,55){\framebox(20,40)}
\put(60,75){\line(1,0){20}}
\put(100,75){\framebox(20,20)}

\put(150,80){$\stackrel{3}{\longrightarrow}$}

\put(175,15){\framebox(140,140)}
\put(265,145){\circle*{4}}
\thicklines
\put(265,145){\line(1,0){50}}
\put(265,145){\line(0,-1){130}}
\put(225,125){\circle*{4}}
\put(225,125){\line(1,0){90}}
\put(225,125){\line(0,-1){110}}
\put(245,105){\circle*{4}}
\put(245,105){\line(1,0){70}}
\put(245,105){\line(0,-1){90}}
\put(305,85){\circle*{4}}
\put(305,85){\line(1,0){10}}
\put(305,85){\line(0,-1){70}}

\put(185,65){\circle*{4}}
\put(185,65){\line(1,0){130}}
\put(185,65){\line(0,-1){50}}

\put(205,45){\circle*{4}}
\put(205,45){\line(1,0){110}}
\put(205,45){\line(0,-1){30}}

\put(285,25){\circle*{4}}
\put(285,25){\line(1,0){30}}
\put(285,25){\line(0,-1){10}}

\thinlines

\put(175,75){\framebox(40,80)}
\put(175,95){\line(1,0){40}}
\put(175,115){\line(1,0){40}}
\put(175,135){\line(1,0){60}}
\put(235,135){\line(0,1){20}}
\put(195,75){\line(0,1){80}}
\put(275,75){\framebox(20,20)}
\put(235,135){\line(1,0){20}}
\put(255,135){\line(0,1){20}}

\put(86,230){\rm X}
\put(86,250){\rm X}
\put(241,250){\rm X}
\put(241,310){\rm X}
\put(0,185){\framebox(140,140)}
\put(70,315){\circle*{4}}
\thicklines
\put(70,315){\line(1,0){70}}
\put(70,315){\line(0,-1){130}}
\put(50,295){\circle*{4}}
\put(50,295){\line(1,0){90}}
\put(50,295){\line(0,-1){110}}
\put(10,275){\circle*{4}}
\put(10,275){\line(1,0){130}}
\put(10,275){\line(0,-1){90}}
\put(130,255){\circle*{4}}
\put(130,255){\line(1,0){10}}
\put(130,255){\line(0,-1){70}}
\put(110,235){\circle*{4}}
\put(110,235){\line(1,0){30}}
\put(110,235){\line(0,-1){50}}
\put(30,215){\circle*{4}}
\put(30,215){\line(1,0){110}}
\put(30,215){\line(0,-1){30}}

\put(90,195){\circle*{4}}
\put(90,195){\line(1,0){50}}
\put(90,195){\line(0,-1){10}}

\thinlines

\put(0,305){\framebox(60,20)}
\put(20,305){\line(0,1){20}}
\put(40,305){\line(0,1){20}}
\put(0,285){\line(1,0){40}}
\put(40,285){\line(0,1){20}}
\put(20,285){\line(0,1){20}}
\put(20,225){\framebox(20,40)}
\put(20,245){\line(1,0){20}}
\put(80,225){\line(0,1){40}}
\put(80,225){\line(1,0){20}}
\put(100,225){\line(0,1){40}}
\put(80,265){\line(1,0){40}}
\put(80,245){\line(1,0){40}}
\put(120,245){\line(0,1){20}}

\put(150,245){$\stackrel{1}{\longrightarrow}$}

\put(175,185){\framebox(140,140)}
\put(265,315){\circle*{4}}
\thicklines
\put(265,315){\line(1,0){50}}
\put(265,315){\line(0,-1){130}}
\put(225,295){\circle*{4}}
\put(225,295){\line(1,0){90}}
\put(225,295){\line(0,-1){110}}
\put(185,275){\circle*{4}}
\put(185,275){\line(1,0){130}}
\put(185,275){\line(0,-1){90}}
\put(305,255){\circle*{4}}
\put(305,255){\line(1,0){10}}
\put(305,255){\line(0,-1){70}}
\put(245,235){\circle*{4}}
\put(245,235){\line(1,0){70}}
\put(245,235){\line(0,-1){50}}
\put(205,215){\circle*{4}}
\put(205,215){\line(1,0){110}}
\put(205,215){\line(0,-1){30}}

\put(285,195){\circle*{4}}
\put(285,195){\line(1,0){30}}
\put(285,195){\line(0,-1){10}}

\thinlines

\put(175,305){\framebox(80,20)}
\put(195,305){\line(0,1){20}}
\put(215,305){\line(0,1){20}}
\put(235,305){\line(0,1){20}}
\put(175,285){\line(1,0){40}}
\put(215,285){\line(0,1){20}}
\put(195,285){\line(0,1){20}}
\put(195,225){\framebox(20,40)}
\put(195,245){\line(1,0){20}}
\put(235,245){\framebox(20,20)}
\put(275,245){\framebox(20,20)}
\end{picture}
\]
\setlength{\unitlength}{.33mm}
\end{Example}


        For any $\beta\in S_n$ and any positive integer $t$, 
we define a rooted, labeled tree
${\mathcal T}_{t}(\beta)$ whose vertices are either labeled by
$\emptyset$ or by a permutation (repetitions allowed). The 
root is labeled by $\beta$. 
If a vertex is labeled by a permutation that has its last descent
weakly smaller than $t$, or is labeled by $\emptyset$, then declare that
vertex to be a leaf. Otherwise, the children of
a vertex are indexed by the output of all ways of marching from that vertex. 
One can check easily that in finitely many steps, this growth 
process terminates,
giving ${\mathcal T}_{t}(\beta)$.
Note that ${\mathcal T}_{t}(\beta)$ 
is a pruning of ${\mathcal T}_{s}(\beta)$ for $t\leq s$. 
Define $K{\mathcal T}_{t}(\beta)$ similarly, 
using instead the $K$-marching operation
(and similarly, 
$K{\mathcal T}_{t}(\beta) \subseteq K{\mathcal T}_{s}(\beta)$).
Finally, if a leaf vertex $v$ is labeled $\pi$, call it a {\dfn $\pi$-leaf}.

We will be particularly interested in the cases that 
$K{\mathcal T}_{s}(\beta)$ or 
${\mathcal T}_{s}(\beta)$ has exactly one labeled leaf 
(i.e. not by $\emptyset$).
The best-behaved cases are when $\beta$ is ``2143-avoiding,'' also known as
``vexillary,'' in which case
$K{\mathcal T}_{s}(\beta)$ has at most one leaf for any $s$.
In the other direction, 
if $\pi$ has a unique descent $\pi(i)>\pi(i+1)$, at $i=s$
(called a {\dfn Grassmannian permutation}), $N\in \naturals$,
and $\beta$ is the {\dfn $N$-stabilization} of of $\pi$, meaning
$$ \beta(i) = i \quad \hbox{ for } i\leq N, \qquad 
        \beta(i) = N+\pi(i-N) \quad \hbox{ for }  i>N, $$
then $K{\mathcal T}_{s}(\beta)$ will have only one labeled leaf, and it
will be labeled $\pi$.

We are now ready to introduce the family of Schubert problems covered
by our main theorem. 
There is a standard operation $\star_n$
on two permutations $\sigma,\alpha\in S_n$. Let
$\sigma\star_{n}\alpha$ be the permutation in $S_{2n}$ whose 
matrix is the direct sum of $G(\sigma)$ and $G(\alpha)$.
For example, $id\star_n\alpha$ is just the $n$-stabilization of $\alpha$.
Let $\sigma\in S_n$ be a permutation whose last descent is $l$, and
let $l\leq t\leq 2n$ be an integer. Suppose that $\alpha\in S_n$ is such that
$K{\mathcal T}_{t}(id\star_{n}\alpha)$ contains a single leaf $v$ with
${\rm label}(v)\neq \emptyset$; let ${\rm label}(v)=\rho$. 
Under these circumstances, call 
$(\sigma,\rho,\pi)\in S_{n}^2 \times S_{\infty}$ a 
{\dfn truncation Schubert problem subjugate to $(t,\alpha)$}.
\begin{Theorem}
\label{thm:main}
If $(\sigma,\rho,\pi)\in S_{n}^2 \times S_{\infty}$ is a 
truncation Schubert problem subjugate to $(t,\alpha)$ then
\[(-1)^{\ell(\sigma)+\ell(\rho)-\ell(\pi)}C_{\sigma,\rho}^{\pi}= 
\mbox{ the number of $\pi$-leaves of 
$K{\mathcal T}_{t}(\sigma\star_n\alpha)$.}\] 
In the cohomology case, we can also say
\[C_{\sigma,\rho}^{\pi}= 
\mbox{ the number of $\pi$-leaves of 
${\mathcal T}_{t}(\sigma\star_n\alpha)$. }\]
\end{Theorem}

\begin{Example}
Let $\sigma=3412$ and $\alpha=3214$ be permutations in $S_4$, 
so $\sigma\star_4\alpha=34127658\in S_8$.
One can check that $K{\mathcal T}_{4}(id\star_4 \alpha)$ 
has a single labeled leaf, labeled by the permutation $12463578$. 
Now $K{\mathcal T}_{4}(\sigma\star_4 \alpha)$ is given in Figure~1 below,
and so by Theorem~\ref{thm:main}:
\begin{figure}
\begin{picture}(420,380)

\put(180,0){$\pi=34127658$}
\put(202,67){\line(0,1){30}}
\put(204,80){{\small $4$}}
\put(185,13){\framebox(48,48)}
\put(185,49){\framebox(12,12)}
\put(191,49){\line(0,1){12}}
\put(185,55){\line(1,0){12}}
\put(209,25){\framebox(6,12)}
\put(209,31){\framebox(12,6)}
\thicklines
\put(230,16){\line(1,0){3}}
\put(230,16){\line(0,-1){3}}
\put(224,34){\line(1,0){9}}
\put(224,34){\line(0,-1){21}}
\put(218,28){\line(1,0){15}}
\put(218,28){\line(0,-1){15}}
\put(212,22){\line(1,0){21}}
\put(212,22){\line(0,-1){9}}
\put(188,46){\line(1,0){45}}
\put(188,46){\line(0,-1){33}}
\put(194,40){\line(1,0){39}}
\put(194,40){\line(0,-1){27}}
\put(200,58){\line(1,0){33}}
\put(200,58){\line(0,-1){45}}
\put(206,52){\line(1,0){27}}
\put(206,52){\line(0,-1){39}}
\thinlines
\put(202,67){\line(-4,1){120}}
\put(202,67){\line(4,1){120}}
\put(168,80){\small $2$}
\put(300,80){\small $2,4$}
\put(50,105){$35127468$}
\put(180,105){$34157268$}
\put(310,105){$35147268$}
\put(50,118){\framebox(48,48)}
\put(50,154){\framebox(12,12)}
\put(50,160){\line(1,0){12}}
\put(56,154){\line(0,1){12}}
\put(68,154){\framebox(6,6)}
\put(68,136){\framebox(6,6)}
\put(80,136){\framebox(6,6)}
\thicklines
\put(65,163){\line(1,0){33}}
\put(65,163){\line(0,-1){45}}
\put(77,157){\line(1,0){21}}
\put(77,157){\line(0,-1){39}}
\put(53,151){\line(1,0){45}}
\put(53,151){\line(0,-1){33}}
\put(59,145){\line(1,0){39}}
\put(59,145){\line(0,-1){27}}
\put(89,139){\line(1,0){9}}
\put(89,139){\line(0,-1){21}}
\put(71,133){\line(1,0){27}}
\put(71,133){\line(0,-1){15}}
\put(83,127){\line(1,0){15}}
\put(83,127){\line(0,-1){9}}
\put(95,121){\line(1,0){3}}
\put(95,121){\line(0,-1){3}}

\thinlines
\put(180,118){\framebox(48,48)}
\put(180,154){\framebox(12,12)}
\put(180,160){\line(1,0){12}}
\put(186,154){\line(0,1){12}}

\put(186,136){\framebox(6,12)}
\put(186,142){\line(1,0){6}}

\put(210,136){\framebox(6,6)}
\thicklines
\put(195,163){\line(1,0){33}}
\put(195,163){\line(0,-1){45}}
\put(201,157){\line(1,0){27}}
\put(201,157){\line(0,-1){39}}
\put(183,151){\line(1,0){45}}
\put(183,151){\line(0,-1){33}}
\put(207,145){\line(1,0){21}}
\put(207,145){\line(0,-1){27}}
\put(219,139){\line(1,0){9}}
\put(219,139){\line(0,-1){21}}
\put(189,133){\line(1,0){39}}
\put(189,133){\line(0,-1){15}}
\put(213,127){\line(1,0){15}}
\put(213,127){\line(0,-1){9}}
\put(225,121){\line(1,0){3}}
\put(225,121){\line(0,-1){3}}
\thinlines

\put(310,118){\framebox(48,48)}
\put(310,154){\framebox(12,12)}
\put(310,160){\line(1,0){12}}
\put(316,154){\line(0,1){12}}

\put(316,136){\framebox(6,12)}
\put(316,142){\line(1,0){6}}

\put(340,136){\framebox(6,6)}
\put(328,154){\framebox(6,6)}
\thicklines
\put(325,163){\line(1,0){33}}
\put(325,163){\line(0,-1){45}}
\put(337,157){\line(1,0){21}}
\put(337,157){\line(0,-1){39}}
\put(313,151){\line(1,0){45}}
\put(313,151){\line(0,-1){33}}
\put(331,145){\line(1,0){27}}
\put(331,145){\line(0,-1){27}}
\put(349,139){\line(1,0){9}}
\put(349,139){\line(0,-1){21}}
\put(319,133){\line(1,0){39}}
\put(319,133){\line(0,-1){15}}
\put(343,127){\line(1,0){15}}
\put(343,127){\line(0,-1){9}}
\put(355,121){\line(1,0){3}}
\put(355,121){\line(0,-1){3}}
\thinlines


\put(72,173){\line(-3,2){50}}
\put(72,173){\line(0,1){30}}
\put(72,173){\line(3,2){50}}
\put(54,185){\small $2$}
\put(74,185){\small $4$}
\put(106,185){\small $2,4$}
\put(0,210){$36125478$}
\put(60,210){$35162478$}
\put(120,210){$36152478$}

\put(60,223){\framebox(48,48)}
\put(60,259){\framebox(12,12)}
\put(60,265){\line(1,0){12}}
\put(66,259){\line(0,1){12}}
\put(78,259){\framebox(6,6)}
\put(66,247){\framebox(6,6)}
\put(78,247){\framebox(6,6)}
\thicklines
\put(75,268){\line(1,0){33}}
\put(75,268){\line(0,-1){45}}
\put(87,262){\line(1,0){21}}
\put(87,262){\line(0,-1){39}}
\put(63,256){\line(1,0){45}}
\put(63,256){\line(0,-1){33}}
\put(93,250){\line(1,0){16}}
\put(93,250){\line(0,-1){27}}
\put(69,244){\line(1,0){39}}
\put(69,244){\line(0,-1){21}}
\put(81,238){\line(1,0){27}}
\put(81,238){\line(0,-1){15}}
\put(99,232){\line(1,0){9}}
\put(99,232){\line(0,-1){9}}
\put(105,226){\line(1,0){3}}
\put(105,226){\line(0,-1){3}}
\thinlines

\put(120,223){\framebox(48,48)}
\put(120,259){\framebox(12,12)}
\put(120,265){\line(1,0){12}}
\put(126,259){\line(0,1){12}}
\put(138,259){\framebox(12,6)}
\put(144,259){\line(0,1){6}}
\put(126,247){\framebox(6,6)}
\put(138,247){\framebox(6,6)}
\thicklines
\put(135,268){\line(1,0){33}}
\put(135,268){\line(0,-1){45}}
\put(153,262){\line(1,0){15}}
\put(153,262){\line(0,-1){39}}
\put(123,256){\line(1,0){45}}
\put(123,256){\line(0,-1){33}}
\put(147,250){\line(1,0){22}}
\put(147,250){\line(0,-1){27}}
\put(129,244){\line(1,0){39}}
\put(129,244){\line(0,-1){21}}
\put(141,238){\line(1,0){27}}
\put(141,238){\line(0,-1){15}}
\put(159,232){\line(1,0){9}}
\put(159,232){\line(0,-1){9}}
\put(165,226){\line(1,0){3}}
\put(165,226){\line(0,-1){3}}
\thinlines

\put(0,223){\framebox(48,48)}
\put(0,259){\framebox(12,12)}
\put(0,265){\line(1,0){12}}
\put(6,259){\line(0,1){12}}
\put(18,259){\framebox(12,6)}
\put(24,259){\line(0,1){6}}
\put(18,241){\framebox(6,6)}
\thicklines
\put(15,268){\line(1,0){33}}
\put(15,268){\line(0,-1){45}}
\put(33,262){\line(1,0){15}}
\put(33,262){\line(0,-1){39}}
\put(3,256){\line(1,0){45}}
\put(3,256){\line(0,-1){33}}

\put(9,250){\line(1,0){39}}
\put(9,250){\line(0,-1){27}}

\put(27,244){\line(1,0){21}}
\put(27,244){\line(0,-1){21}}
\put(21,238){\line(1,0){27}}
\put(21,238){\line(0,-1){15}}
\put(39,232){\line(1,0){9}}
\put(39,232){\line(0,-1){9}}
\put(45,226){\line(1,0){3}}
\put(45,226){\line(0,-1){3}}
\thinlines

\put(202,173){\line(0,1){30}}
\put(204,185){\small $4$}
\put(180,210){$34165278$}

\put(180,223){\framebox(48,48)}
\put(180,259){\framebox(12,12)}
\put(180,265){\line(1,0){12}}
\put(186,259){\line(0,1){12}}

\put(186,241){\framebox(6,12)}
\put(186,247){\line(1,0){6}}

\put(204,247){\framebox(6,6)}
\thicklines
\put(195,268){\line(1,0){33}}
\put(195,268){\line(0,-1){45}}
\put(201,262){\line(1,0){27}}
\put(201,262){\line(0,-1){39}}
\put(183,256){\line(1,0){45}}
\put(183,256){\line(0,-1){33}}
\put(213,250){\line(1,0){15}}
\put(213,250){\line(0,-1){27}}
\put(207,244){\line(1,0){21}}
\put(207,244){\line(0,-1){21}}
\put(189,238){\line(1,0){39}}
\put(189,238){\line(0,-1){15}}
\put(219,232){\line(1,0){9}}
\put(219,232){\line(0,-1){9}}
\put(225,226){\line(1,0){3}}
\put(225,226){\line(0,-1){3}}
\thinlines

\put(332,173){\line(0,1){30}}
\put(332,173){\line(-3,2){50}}
\put(332,173){\line(3,2){50}}
\put(314,185){\small $2$}
\put(334,185){\small $4$}
\put(364,185){\small $2,4$}
\put(240,210){$36145278$}
\put(310,210){$35164278$}
\put(370,210){$36154278$}

\put(310,223){\framebox(48,48)}
\put(310,259){\framebox(12,12)}
\put(310,265){\line(1,0){12}}
\put(316,259){\line(0,1){12}}

\put(316,241){\framebox(6,12)}
\put(316,247){\line(1,0){6}}

\put(328,247){\framebox(6,6)}
\put(328,259){\framebox(6,6)}
\thicklines
\put(325,268){\line(1,0){33}}
\put(325,268){\line(0,-1){45}}
\put(337,262){\line(1,0){21}}
\put(337,262){\line(0,-1){39}}
\put(313,256){\line(1,0){45}}
\put(313,256){\line(0,-1){33}}
\put(343,250){\line(1,0){15}}
\put(343,250){\line(0,-1){27}}
\put(331,244){\line(1,0){27}}
\put(331,244){\line(0,-1){21}}
\put(319,238){\line(1,0){39}}
\put(319,238){\line(0,-1){15}}
\put(349,232){\line(1,0){9}}
\put(349,232){\line(0,-1){9}}
\put(355,226){\line(1,0){3}}
\put(355,226){\line(0,-1){3}}
\thinlines

\put(370,223){\framebox(48,48)}
\put(370,259){\framebox(12,12)}
\put(370,265){\line(1,0){12}}
\put(376,259){\line(0,1){12}}

\put(376,241){\framebox(6,12)}
\put(376,247){\line(1,0){6}}

\put(388,247){\framebox(6,6)}
\put(388,259){\framebox(12,6)}
\put(394,259){\line(0,1){6}}
\thicklines
\put(385,268){\line(1,0){33}}
\put(385,268){\line(0,-1){45}}
\put(403,262){\line(1,0){15}}
\put(403,262){\line(0,-1){39}}
\put(373,256){\line(1,0){45}}
\put(373,256){\line(0,-1){33}}
\put(397,250){\line(1,0){21}}
\put(397,250){\line(0,-1){27}}
\put(391,244){\line(1,0){27}}
\put(391,244){\line(0,-1){21}}
\put(379,238){\line(1,0){39}}
\put(379,238){\line(0,-1){15}}
\put(409,232){\line(1,0){9}}
\put(409,232){\line(0,-1){9}}
\put(415,226){\line(1,0){3}}
\put(415,226){\line(0,-1){3}}
\thinlines

\put(240,223){\framebox(48,48)}
\put(240,259){\framebox(12,12)}
\put(240,265){\line(1,0){12}}
\put(246,259){\line(0,1){12}}

\put(246,241){\framebox(6,12)}
\put(246,247){\line(1,0){6}}

\put(258,259){\framebox(12,6)}
\put(264,259){\line(0,1){6}}
\thicklines
\put(255,268){\line(1,0){33}}
\put(255,268){\line(0,-1){45}}
\put(273,262){\line(1,0){15}}
\put(273,262){\line(0,-1){39}}
\put(243,256){\line(1,0){45}}
\put(243,256){\line(0,-1){33}}
\put(261,250){\line(1,0){27}}
\put(261,250){\line(0,-1){27}}
\put(267,244){\line(1,0){21}}
\put(267,244){\line(0,-1){21}}
\put(249,238){\line(1,0){39}}
\put(249,238){\line(0,-1){15}}
\put(279,232){\line(1,0){9}}
\put(279,232){\line(0,-1){9}}
\put(285,226){\line(1,0){3}}
\put(285,226){\line(0,-1){3}}
\thinlines


\put(24,278){\line(0,1){30}}
\put(26,292){$1$}
\put(24,278){\line(2,1){60}}
\put(67,292){$4$}
\put(24,278){\line(4,1){120}}
\put(118,292){$1,4$}
\put(0,312){$46123578$}
\put(60,312){$36142578$}
\put(120,312){$46132578$}

\put(0,323){\framebox(48,48)}
\put(0,359){\framebox(18,12)}
\put(0,365){\line(1,0){18}}
\put(6,359){\line(0,1){12}}
\put(12,359){\line(0,1){12}}
\put(24,359){\framebox(6,6)}

\thicklines
\put(21,368){\line(1,0){27}}
\put(21,368){\line(0,-1){45}}
\put(33,362){\line(1,0){15}}
\put(33,362){\line(0,-1){39}}
\put(3,356){\line(1,0){45}}
\put(3,356){\line(0,-1){33}}

\put(9,350){\line(1,0){39}}
\put(9,350){\line(0,-1){27}}

\put(15,344){\line(1,0){33}}
\put(15,344){\line(0,-1){21}}
\put(27,338){\line(1,0){21}}
\put(27,338){\line(0,-1){15}}
\put(39,332){\line(1,0){9}}
\put(39,332){\line(0,-1){9}}
\put(45,326){\line(1,0){3}}
\put(45,326){\line(0,-1){3}}
\thinlines


\put(120,323){\framebox(48,48)}
\put(120,359){\framebox(18,12)}
\put(120,365){\line(1,0){18}}
\put(126,359){\line(0,1){12}}
\put(132,359){\line(0,1){12}}
\put(144,359){\framebox(6,6)}
\put(126,347){\framebox(6,6)}

\thicklines
\put(141,368){\line(1,0){27}}
\put(141,368){\line(0,-1){45}}
\put(153,362){\line(1,0){15}}
\put(153,362){\line(0,-1){39}}
\put(123,356){\line(1,0){45}}
\put(123,356){\line(0,-1){33}}

\put(135,350){\line(1,0){33}}
\put(135,350){\line(0,-1){27}}

\put(129,344){\line(1,0){39}}
\put(129,344){\line(0,-1){21}}
\put(147,338){\line(1,0){21}}
\put(147,338){\line(0,-1){15}}
\put(159,332){\line(1,0){9}}
\put(159,332){\line(0,-1){9}}
\put(165,326){\line(1,0){3}}
\put(165,326){\line(0,-1){3}}
\thinlines

\put(60,323){\framebox(48,48)}
\put(60,359){\framebox(12,12)}
\put(60,365){\line(1,0){12}}
\put(66,359){\line(0,1){12}}
\put(78,359){\framebox(12,6)}
\put(84,359){\line(0,1){6}}
\put(66,347){\framebox(6,6)}

\thicklines
\put(75,368){\line(1,0){33}}
\put(75,368){\line(0,-1){45}}
\put(93,362){\line(1,0){15}}
\put(93,362){\line(0,-1){39}}
\put(63,356){\line(1,0){45}}
\put(63,356){\line(0,-1){33}}

\put(81,350){\line(1,0){27}}
\put(81,350){\line(0,-1){27}}

\put(69,344){\line(1,0){39}}
\put(69,344){\line(0,-1){21}}
\put(87,338){\line(1,0){21}}
\put(87,338){\line(0,-1){15}}
\put(99,332){\line(1,0){9}}
\put(99,332){\line(0,-1){9}}
\put(105,326){\line(1,0){3}}
\put(105,326){\line(0,-1){3}}
\thinlines

\put(202,278){\line(0,1){30}}
\put(204,292){$3$}
\put(180,312){$34261578$}

\put(180,323){\framebox(48,48)}
\put(180,359){\framebox(12,12)}
\put(180,365){\line(1,0){12}}
\put(186,359){\line(0,1){12}}

\put(180,347){\framebox(6,12)}
\put(180,353){\line(1,0){6}}

\put(204,347){\framebox(6,6)}
\thicklines
\put(195,368){\line(1,0){33}}
\put(195,368){\line(0,-1){45}}
\put(201,362){\line(1,0){27}}
\put(201,362){\line(0,-1){39}}
\put(189,356){\line(1,0){39}}
\put(189,356){\line(0,-1){33}}
\put(213,350){\line(1,0){15}}
\put(213,350){\line(0,-1){27}}
\put(183,344){\line(1,0){45}}
\put(183,344){\line(0,-1){21}}
\put(207,338){\line(1,0){21}}
\put(207,338){\line(0,-1){15}}
\put(219,332){\line(1,0){9}}
\put(219,332){\line(0,-1){9}}
\put(225,326){\line(1,0){3}}
\put(225,326){\line(0,-1){3}}
\thinlines

\put(262,278){\line(0,1){30}}
\put(264,292){$3$}
\put(240,312){$36241578$}

\put(240,323){\framebox(48,48)}
\put(240,359){\framebox(12,12)}
\put(240,365){\line(1,0){12}}
\put(246,359){\line(0,1){12}}

\put(240,347){\framebox(6,12)}
\put(240,353){\line(1,0){6}}

\put(258,359){\framebox(12,6)}
\put(264,359){\line(0,1){6}}
\thicklines
\put(255,368){\line(1,0){33}}
\put(255,368){\line(0,-1){45}}
\put(273,362){\line(1,0){15}}
\put(273,362){\line(0,-1){39}}
\put(249,356){\line(1,0){39}}
\put(249,356){\line(0,-1){33}}
\put(262,350){\line(1,0){27}}
\put(262,350){\line(0,-1){27}}
\put(243,344){\line(1,0){45}}
\put(243,344){\line(0,-1){21}}
\put(267,338){\line(1,0){21}}
\put(267,338){\line(0,-1){15}}
\put(279,332){\line(1,0){9}}
\put(279,332){\line(0,-1){9}}
\put(285,326){\line(1,0){3}}
\put(285,326){\line(0,-1){3}}
\thinlines

\put(334,278){\line(0,1){30}}
\put(336,292){$3$}
\put(310,312){$35261478$}

\put(310,323){\framebox(48,48)}
\put(310,359){\framebox(12,12)}
\put(310,365){\line(1,0){12}}
\put(316,359){\line(0,1){12}}

\put(310,347){\framebox(6,12)}
\put(310,353){\line(1,0){6}}

\put(328,359){\framebox(6,6)}

\put(328,347){\framebox(6,6)}
\thicklines
\put(325,368){\line(1,0){33}}
\put(325,368){\line(0,-1){45}}
\put(337,362){\line(1,0){21}}
\put(337,362){\line(0,-1){39}}
\put(319,356){\line(1,0){39}}
\put(319,356){\line(0,-1){33}}
\put(344,350){\line(1,0){15}}
\put(344,350){\line(0,-1){27}}
\put(313,344){\line(1,0){45}}
\put(313,344){\line(0,-1){21}}
\put(331,338){\line(1,0){27}}
\put(331,338){\line(0,-1){15}}
\put(349,332){\line(1,0){9}}
\put(349,332){\line(0,-1){9}}
\put(355,326){\line(1,0){3}}
\put(355,326){\line(0,-1){3}}
\thinlines

\put(394,278){\line(0,1){30}}
\put(396,292){$3$}
\put(370,312){$36251478$}

\put(370,323){\framebox(48,48)}
\put(370,359){\framebox(12,12)}
\put(370,365){\line(1,0){12}}
\put(376,359){\line(0,1){12}}

\put(370,347){\framebox(6,12)}
\put(370,353){\line(1,0){6}}

\put(388,359){\framebox(12,6)}
\put(394,359){\line(0,1){6}}
\put(388,347){\framebox(6,6)}
\thicklines
\put(385,368){\line(1,0){33}}
\put(385,368){\line(0,-1){45}}
\put(403,362){\line(1,0){15}}
\put(403,362){\line(0,-1){39}}
\put(379,356){\line(1,0){39}}
\put(379,356){\line(0,-1){33}}
\put(398,350){\line(1,0){21}}
\put(398,350){\line(0,-1){27}}
\put(373,344){\line(1,0){45}}
\put(373,344){\line(0,-1){21}}
\put(391,338){\line(1,0){27}}
\put(391,338){\line(0,-1){15}}
\put(409,332){\line(1,0){9}}
\put(409,332){\line(0,-1){9}}
\put(415,326){\line(1,0){3}}
\put(415,326){\line(0,-1){3}}
\thinlines
\end{picture}
\caption{
The tree $K{\mathcal T}_{4}(34127657)$; see Example~3.}
\end{figure}
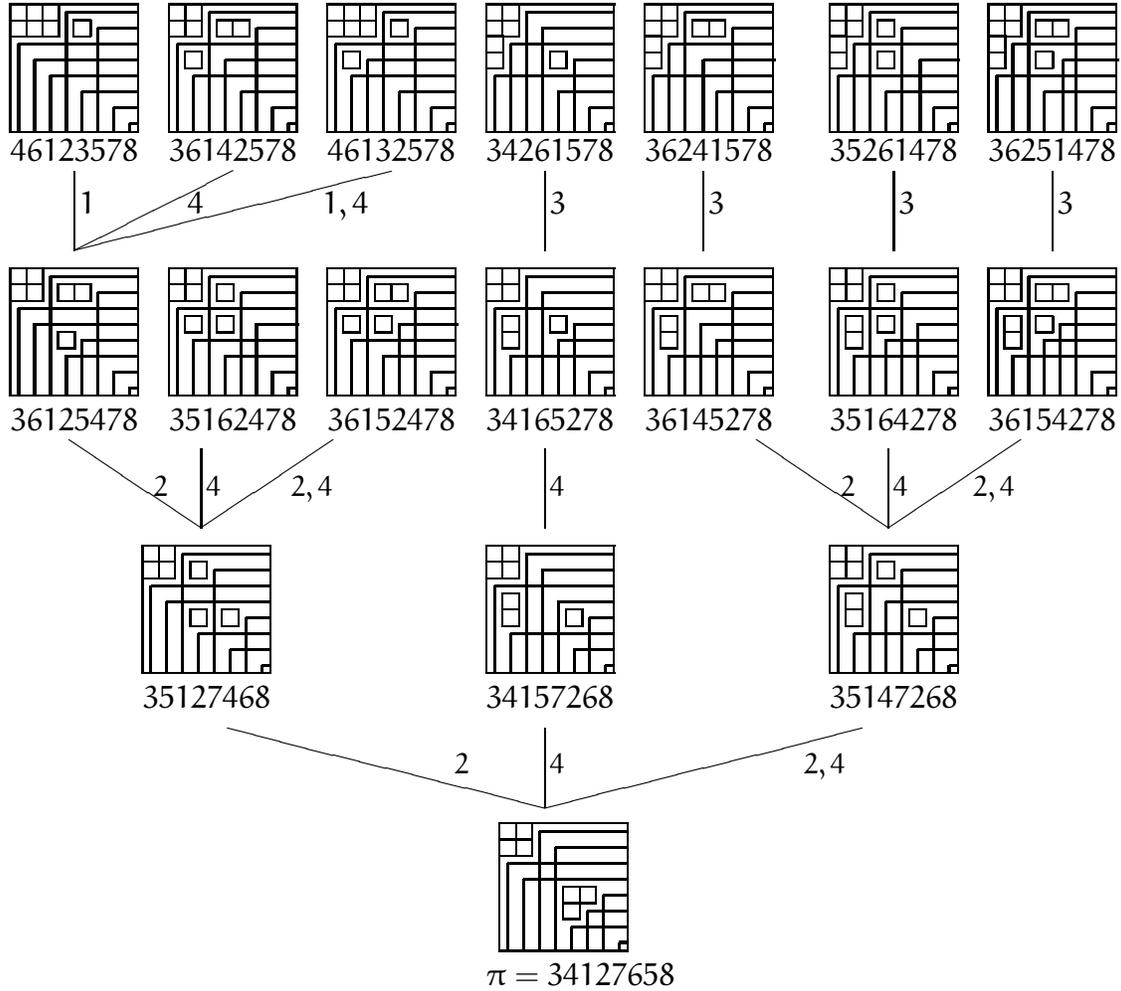
\begin{multline}\nonumber
[{\mathcal O}_{X_{3412}}]\cdot[{\mathcal O}_{X_{12463578}}] =  
[{\mathcal O}_{X_{46123578}}]+[{\mathcal O}_{X_{36142578}}]
+[{\mathcal O}_{X_{35162478}}]
+[{\mathcal O}_{X_{34261578}}]\\ \nonumber
-[{\mathcal O}_{X_{46132578}}]
-[{\mathcal O}_{X_{36152478}}]
-[{\mathcal O}_{X_{36241578}}]
-[{\mathcal O}_{X_{35261478}}]
+[{\mathcal O}_{X_{36251478}}]
\end{multline}
where the expansion (\ref{eqn:main_def}) 
has been done in the case ${\rm Flags}({\mathbb C}^8)$.
\end{Example}

As mentioned before the theorem, 
and spelled out in the corollary below,
one family of truncation Schubert problems comes from
Grassmannian permutations. 
In the cohomology case, 
these were given a (different) positive 
combinatorial formula by M.~Kogan~\cite{kogan:schur_schub}.

\begin{Corollary}
\label{cor:kogan}
Let $\sigma\in S_n$ have last descent at $l$ and 
let $\rho$ be a Grassmannian permutation with unique descent
at $t$, where $l \leq t\leq n$.
Then for any $\pi\in S_{\infty}$,
\[(-1)^{\ell(\sigma)+\ell(\rho)-\ell(\pi)}C_{\sigma,\rho}^{\pi}= 
\mbox{ the number of $\pi$-leaves of 
$K{\mathcal T}_{t}(\sigma\!\star_n\!\rho)$. }\]
In the cohomology case (treated in \cite{kogan:schur_schub}),
\[C_{\sigma,\rho}^{\pi}= 
\mbox{ the number of $\pi$-leaves of 
${\mathcal T}_{t}(\sigma\!\star_n\!\rho)$}. \]
\end{Corollary}

	If we also assume that
$\sigma$ is grassmannian and moreover $t=l$,
the first conclusion of Corollary~\ref{cor:kogan} computes 
the $K$-theory Littlewood-Richardson coefficients of~\cite{buch:KLR},
while the second conclusion computes the classical Littlewood-Richardson
coefficients.

\begin{Example}
Let $\sigma=321\in S_3, \rho=132$ and $\sigma\star_3\rho=321465$.
The tree $K{\mathcal T}_{2}(321465)$ is given in Figure~2 below.
Using this, the expansion (\ref{eqn:main_def}) for ${\rm Flags}({\mathbb C}^6)$
is 
\[[{\mathcal O}_{X_{321}}]\cdot[{\mathcal O}_{X_{132}}]
=[{\mathcal O}_{X_{421356}}]+[{\mathcal O}_{X_{341256}}]-[{\mathcal O}_{X_{431256}}]
.\]
\begin{figure}
\begin{picture}(420,300)

\put(205,75){{\small $4$}}
\put(175,0){$\pi=321465$}
\put(185,15){\framebox(36,36)}
\put(185,39){\framebox(6,12)}
\put(185,45){\framebox(12,6)}
\put(209,21){\framebox(6,6)}
\thicklines
\put(200,48){\line(1,0){21}}
\put(200,48){\line(0,-1){33}}
\put(194,42){\line(1,0){27}}
\put(194,42){\line(0,-1){27}}
\put(188,36){\line(1,0){33}}
\put(188,36){\line(0,-1){21}}
\put(206,30){\line(1,0){15}}
\put(206,30){\line(0,-1){15}}
\put(218,24){\line(1,0){3}}
\put(218,24){\line(0,-1){9}}
\put(212,18){\line(1,0){9}}
\put(212,18){\line(0,-1){3}}
\thinlines
\put(200,62){\line(0,1){30}}

\put(185,100){$321546$}
\put(185,115){\framebox(36,36)}
\put(185,139){\framebox(6,12)}
\put(185,145){\framebox(12,6)}
\put(203,127){\framebox(6,6)}
\thicklines
\put(200,148){\line(1,0){21}}
\put(200,148){\line(0,-1){33}}
\put(194,142){\line(1,0){27}}
\put(194,142){\line(0,-1){27}}
\put(188,136){\line(1,0){33}}
\put(188,136){\line(0,-1){21}}
\put(212,130){\line(1,0){9}}
\put(212,130){\line(0,-1){15}}
\put(206,124){\line(1,0){15}}
\put(206,124){\line(0,-1){9}}
\put(218,118){\line(1,0){3}}
\put(218,118){\line(0,-1){3}}

\thinlines
\put(175,160){\line(-5,1){140}}

\put(175,160){\line(-5,2){75}}
\put(200,160){\line(-3,2){45}}
\put(200,160){\line(0,1){30}}
\put(200,160){\line(3,2){45}}
\put(230,160){\line(5,1){140}}
\put(230,160){\line(5,2){75}}

\put(80,180){{\small $1$}}
\put(130,180){{\small $2$}}
\put(180,180){{\small $3$}}
\put(205,180){{\small $1,2$}}
\put(245,180){{\small $2,3$}}
\put(300,180){{\small $1,3$}}
\put(370,180){{\small $1,2,3$}}
\thinlines

\put(20,195){{$421356$}}

\put(20,210){\framebox(36,36)}
\put(20,234){\framebox(6,12)}
\put(20,240){\framebox(18,6)}
\put(32,240){\line(0,1){6}}
\thicklines
\put(41,243){\line(1,0){15}}
\put(41,243){\line(0,-1){33}}
\put(29,237){\line(1,0){27}}
\put(29,237){\line(0,-1){27}}
\put(23,231){\line(1,0){33}}
\put(23,231){\line(0,-1){21}}

\put(35,225){\line(1,0){21}}
\put(35,225){\line(0,-1){15}}

\put(47,219){\line(1,0){9}}
\put(47,219){\line(0,-1){9}}
\put(53,213){\line(1,0){3}}
\put(53,213){\line(0,-1){3}}
\thinlines

\put(80,195){{$341256$}}
\put(80,210){\framebox(36,36)}
\put(80,234){\framebox(12,12)}
\put(80,240){\line(1,0){12}}
\put(86,234){\line(0,1){12}}
\thicklines
\put(95,243){\line(1,0){21}}
\put(95,243){\line(0,-1){33}}

\put(101,237){\line(1,0){15}}
\put(101,237){\line(0,-1){27}}

\put(83,231){\line(1,0){33}}
\put(83,231){\line(0,-1){21}}

\put(89,225){\line(1,0){27}}
\put(89,225){\line(0,-1){15}}

\put(107,219){\line(1,0){9}}
\put(107,219){\line(0,-1){9}}
\put(113,213){\line(1,0){3}}
\put(113,213){\line(0,-1){3}}
\thinlines


\put(130,195){{$324156$}}
\put(130,210){\framebox(36,36)}
\put(130,228){\framebox(6,18)}
\put(130,234){\line(1,0){6}}
\put(130,240){\framebox(12,6)}

\thicklines
\put(145,243){\line(1,0){21}}
\put(145,243){\line(0,-1){33}}

\put(139,237){\line(1,0){27}}
\put(139,237){\line(0,-1){27}}

\put(151,231){\line(1,0){15}}
\put(151,231){\line(0,-1){21}}

\put(133,225){\line(1,0){33}}
\put(133,225){\line(0,-1){15}}

\put(157,219){\line(1,0){9}}
\put(157,219){\line(0,-1){9}}
\put(163,213){\line(1,0){3}}
\put(163,213){\line(0,-1){3}}
\thinlines


\put(180,195){{$431256$}}

\put(180,210){\framebox(36,36)}
\put(180,234){\framebox(12,12)}
\put(180,240){\framebox(18,6)}
\put(180,234){\framebox(6,12)}
\thicklines
\put(201,243){\line(1,0){15}}
\put(201,243){\line(0,-1){33}}

\put(195,237){\line(1,0){21}}
\put(195,237){\line(0,-1){27}}

\put(183,231){\line(1,0){33}}
\put(183,231){\line(0,-1){21}}

\put(189,225){\line(1,0){27}}
\put(189,225){\line(0,-1){15}}

\put(207,219){\line(1,0){9}}
\put(207,219){\line(0,-1){9}}
\put(213,213){\line(1,0){3}}
\put(213,213){\line(0,-1){3}}
\thinlines

\put(230,195){{$342156$}}
\put(230,210){\framebox(36,36)}
\put(230,234){\framebox(12,12)}
\put(230,240){\framebox(12,6)}
\put(230,228){\framebox(6,18)}

\thicklines
\put(245,243){\line(1,0){21}}
\put(245,243){\line(0,-1){33}}
\put(251,237){\line(1,0){15}}
\put(251,237){\line(0,-1){27}}

\put(239,231){\line(1,0){27}}
\put(239,231){\line(0,-1){21}}

\put(233,225){\line(1,0){33}}
\put(233,225){\line(0,-1){15}}

\put(257,219){\line(1,0){9}}
\put(257,219){\line(0,-1){9}}

\put(263,213){\line(1,0){3}}
\put(263,213){\line(0,-1){3}}
\thinlines

\put(290,195){{$423156$}}
\put(290,210){\framebox(36,36)}
\put(290,228){\framebox(6,18)}
\put(290,240){\framebox(18,6)}
\put(290,234){\line(1,0){6}}
\put(302,240){\line(0,1){6}}
\thicklines
\put(311,243){\line(1,0){15}}
\put(311,243){\line(0,-1){33}}

\put(299,237){\line(1,0){27}}
\put(299,237){\line(0,-1){27}}

\put(305,231){\line(1,0){21}}
\put(305,231){\line(0,-1){21}}

\put(317,219){\line(1,0){9}}
\put(317,219){\line(0,-1){9}}

\put(293,225){\line(1,0){33}}
\put(293,225){\line(0,-1){15}}

\put(323,213){\line(1,0){3}}
\put(323,213){\line(0,-1){3}}
\thinlines

\put(360,195){{$432156$}}
\put(360,210){\framebox(36,36)}
\put(360,234){\framebox(12,12)}
\put(360,240){\framebox(18,6)}
\put(360,228){\framebox(6,18)}

\thicklines
\put(381,243){\line(1,0){15}}
\put(381,243){\line(0,-1){33}}
\put(375,237){\line(1,0){21}}
\put(375,237){\line(0,-1){27}}
\put(369,231){\line(1,0){27}}
\put(369,231){\line(0,-1){21}}
\put(387,219){\line(1,0){9}}
\put(387,219){\line(0,-1){9}}

\put(363,225){\line(1,0){33}}
\put(363,225){\line(0,-1){15}}

\put(393,213){\line(1,0){3}}
\put(393,213){\line(0,-1){3}}
\thinlines

\put(145,250){\line(0,1){30}}
\put(143,285){$\emptyset$}

\put(247,250){\line(0,1){30}}
\put(245,285){$\emptyset$}

\put(305,250){\line(0,1){30}}
\put(303,285){$\emptyset$}

\put(375,250){\line(0,1){30}}
\put(373,285){$\emptyset$}
\end{picture}
\caption{The tree $K{\mathcal T}_{2}(321465)$; see Example~4.}
\end{figure}
\end{Example}
There is an isomorphism of $\Fln$ to itself induced by sending each vector
subspace $V$ to its orthogonal complement $V^{\perp}$ (with respect to
an arbitrarily chosen bilinear form). The induced automorphism
of $K(\Fln)$ gives the symmetry
$C_{\sigma,\rho}^{\pi}=C_{w_0 \sigma w_0, w_0\rho w_0}^{w_0 \pi w_0}$.
This observation, combined with the corollary (or the theorem), gives, e.g., a
subtraction-free formula also for the Schubert numbers
$C_{\sigma,\rho}^{\pi}$ where $\rho$ is Grassmannian of descent $t$
which is weakly {\em smaller}
than the {\em first} descent of $\sigma$.

Theorem~\ref{thm:main} also handles some new (but apparently limited) cases
of Schubert problems $(\sigma,\rho)$ where neither $\sigma$ nor $\rho$
are Grassmannian permutations. This differs from other formulas, see, 
e.g.,~\cite{bergeron.sottile:duke,buch:KLR,buch.sottile.yong:quiver,
kogan:schur_schub,lascoux,lenart,lenart.sottile:pieri,manivel,
sottile:pieri,vakil:LR}. 

\begin{Example}
The tree $K{\mathcal T}_{7}(123459876{\underline{10}})$ has a single 
leaf indexed by a permutation, and that permutation is 
$123469857{\underline{10}}$. Hence a 
product of $[{\mathcal O}_{X_{123469857{\underline{10}}}}]$ with
any $[{\mathcal O}_{X_{\rho}}]$ where $\rho\in S_5$ 
is covered by Theorem~\ref{thm:main}, and in particular in 
$K({\rm Flags}({\mathbb C}^{10}))$: 
\[[{\mathcal O}_{X_{123469857{\underline{10}}}}]\cdot
[{\mathcal O}_{X_{{41352}}}]
=[{\mathcal O}_{X_{{413629857{\underline{10}}}}}]
+[{\mathcal O}_{X_{413569827{\underline{10}}}}]
-[{\mathcal O}_{X_{413659827{\underline{10}}}}]\]
is a nontrivial expansion which is not computed 
any previously known (subtraction-free) multiplication formula. We remark
that Theorem~\ref{thm:main} is the first to give a positive formula for
even the cohomology expansion in $H^{*}({\rm Flags}({\mathbb C}^{10}))$:
\[[X_{123469857{\underline{10}}}]\cdot
[{X_{{41352}}}]
=[{X_{{413629857{\underline{10}}}}}]
+[{X_{413569827{\underline{10}}}}].\]
\end{Example}

\section{Proof of Theorem~1 and Corollary~1}

We begin by recalling A.~Lascoux and M.\!~-P.~Sch\"{u}tzenberger's 
{\dfn Grothendieck polynomials}~\cite{lascoux.schutzenberger:hopf},
albeit via a rather unconventional definition.
Let $X=\{x_1,x_2,\ldots \ \}$ be a collection of
commuting independent variables. 
To each $\pi \in S_\infty$,
there is an associated Grothendieck polynomial in the $\{x_i\}$,
and these polynomials satisfy the following crucial recursion:

\begin{Theorem}
\label{DefThm:trans}
{\rm (cf.~\cite{lascoux,lenart})}\label{lem:Ktransition}
For any permutation $\gamma\in S_{\infty}$ with last descent $g$, 
let $m>g$ be the largest integer such that
$\gamma(m)<\gamma(g)$ and set $\gamma'=\gamma t_{g\leftrightarrow m}$. 
Suppose that $1\leq i_1<i_2<\ldots<i_s<g$ are the positions such that
$\ell(\gamma' t_{i_j\leftrightarrow g})=\ell(\gamma')+1$. Then
the ($K$-theory) transition formula of A.~Lascoux~\cite{lascoux} 
(we give the formulation~\cite[Cor.~3.10]{lenart}) holds:
\begin{equation}
\label{eqn:K-transition}
\groth_{\gamma}(X) = \groth_{\gamma'}(X)
+(x_g -1)\big[\groth_{\gamma'}(X)\cdot ({\rm I}-t_{i_1\leftrightarrow g})\cdots
({\rm I}-t_{i_s\leftrightarrow g})\big],
\end{equation}
where $t_{j\leftrightarrow l}$ acts on the $\{\groth_{\xi}(X)\}$ by 
$\groth_{\xi}(X)\cdot t_{j\leftrightarrow l}
= \groth_{\xi t_{j\leftrightarrow l}}(X)$ 
and ${\rm I}$ acts as the identity operator. 
\end{Theorem}

This, and the base case $\groth_{id}(X) = 1$, 
uniquely determine the Grothendieck polynomials (the usual definition
is via isobaric divided difference operators). Together, these polynomials
form a ${\mathbb Z}$-linear basis of ${\mathbb Z}[X]$ and satisfy
\begin{equation*}
\label{eqn:groth_mult}
\groth_{\sigma}(X)\groth_{\rho}(X)
=\sum_{\pi\in S_{\infty}}C_{\sigma,\rho}^{\pi}\groth_{\pi}(X).
\end{equation*} 

For any positive integer $t$, define the {\dfn truncation homomorphism} 
$r_{t}:{\mathbb Z}[X]\to {\mathbb Z}[X]$ by 
\linebreak $r_{t}(f(X))=f(x_1,\ldots,x_t,0,0,\ldots)$.  

Theorem~\ref{thm:main} is immediate from the second formula of the
following result:

\begin{Theorem}
\label{thm:2}
For any $\gamma\in S_{\infty}$, we have
\[r_{t}(\groth_{\gamma}(X))=\sum_{v} 
(-1)^{\ell(\gamma)-\ell({\rm label}(v))} \groth_{{\rm label}(v)}(X),\]
where the sum is over all leaves $v$ of $K{\mathcal T}_{t}(\gamma)$
such that ${\rm label}(v)\neq \emptyset$.

If $\sigma\in S_n$ has its last descent weakly smaller than $t$, 
and $\alpha\in S_n$ is arbitrary, then 
\[\groth_{\sigma}(X) \, r_{t}(\groth_{id\star_n \alpha}(X))
=\sum_{v}(-1)^{\ell(\sigma)+\ell(\alpha)
-\ell({\rm label}(v))}\groth_{{\rm label}(v)}(X),\] 
where the sum is over all leaves $v$ of 
$K{\mathcal T}_{t}(\sigma\star_n \alpha)$ such that 
${\rm label}(v)\neq \emptyset$.
\end{Theorem}

\begin{proof}
  To expand $r_{t}(\groth_{\gamma}(X))$ we will need the
  following lemma, which connects the diagram moves of Section~2 to
  the ($K$-theory) transition formula in Theorem~\ref{lem:Ktransition}.
  It also gives an alternative form of a substitution formula
  of~\cite{lenart.robinson.sottile}.

\begin{Lemma} 
\label{lemma:lrs_like}
{\rm (cf.~\cite{lenart.robinson.sottile})}
Under the assumptions of Theorem~\ref{DefThm:trans}
we have that
\begin{itemize}
\item[(i)] $(i_1,\gamma(i_1)),\ldots, (i_s,\gamma(i_s))$ are the pivots 
of the maximal corner $(g,\gamma^{-1}(m))\in D(\gamma)$;
\item[(ii)] the following formula holds:
\begin{equation*}
\label{eqn:our_transition}
r_{g-1}(\groth_{\gamma}(X))=
\sum_{\gamma\stackrel{{\mathcal I}}{\longrightarrow}\tau}
(-1)^{\ell(\gamma)-\ell(\tau)}r_{g-1}(\groth_{\tau}(X))
\end{equation*}
where the summation ranges over all subsets 
${\mathcal I}$ of $\{i_1,\ldots,i_s\}$.
\end{itemize}
\end{Lemma}

\begin{proof}
Observe that the diagram of $\gamma'$ differs from the diagram of
$\gamma$ only in that the maximal accessible box
of $\gamma$ has been removed 
(and thus there is a dot of $G(\gamma')$ in that position instead).\footnote{%
This seems to be the main reason to work with diagrams of permutations
rather than inversion sets (which more easily generalize to 
other root systems).}
Now, for any index 
$1\leq a< g$, $\ell(\gamma' t_{a\leftrightarrow g})=\ell(\gamma')+1$ 
holds if and only if $\gamma'(a)<\gamma'(g)$ and the rectangle defined by
$(a,\gamma'(a))$ and $(g,\gamma'(g))$ contains no other dots of
$G(\gamma')$; that is, if and only if $(a,\gamma'(a))$ is a pivot of
the maximal accessible box of $\gamma$.  Hence (i) holds.

Thus in view of (i), conclusion (ii) follows easily by expanding
the $K$-transition formula from Theorem~\ref{lem:Ktransition}, observing that
$(-1)^{\# {\mathcal I}}=(-1)^{\ell(\gamma)-\ell(\tau)}$, and setting $x_g = 0$.
\end{proof}

From the above lemma, we have
\begin{eqnarray}
 \label{eqn:ind_step}
 r_{g-1}(\groth_{\gamma}(X)) & = & 
 \sum_{\gamma\stackrel{I}{\longrightarrow}\tau}
 (-1)^{\ell(\gamma)-\ell(\tau)} r_{g-1}(\groth_{\tau}(X)) \nonumber \\ 
 & = & \sum_{\mathcal A}(-1)^{\ell(\gamma)-\ell(\tau)}\groth_{\tau}
 +\sum_{\mathcal B}(-1)^{\ell(\gamma)-\ell(\tau)}r_{g-1}(\groth_{\tau}(X)),
\end{eqnarray}
where ${\mathcal A}$ consists of those $\tau$ appearing from marching from
$\gamma$ that have last descent at $g-1$ or smaller,
and ${\mathcal B}$ consists of those
that still have their last descent at $g$.
It is not hard to see that a finite number of iterations of marches 
from $\gamma$ results in $\emptyset$ or a permutation with
last descent weakly smaller than $g-1$. Thus after repeated application
of Lemma~\ref{lemma:lrs_like} on (\ref{eqn:ind_step}) we expand 
$r_{g-1}(\groth_{\gamma}(X))$ into the sum of Grothendieck polynomials
indexed by such permutations (in particular we have just used the fact that
if a permutation $\tau$ has last descent $t$ and has no pivots, then
$r_{t-1}(\groth_{\tau}(X))=0$). 
Therefore, the first conclusion follows by iterating the
operation of setting the variables $x_{g-1},x_{g-2},\ldots,x_{t+1}$ 
to zero in succession. 

Only a little more is necessary for the second conclusion of the theorem.
Using Theorem~\ref{DefThm:trans} and induction,
it is easy to check that
\[
\groth_{\sigma}(X)\groth_{id \star_n \rho}(X)=\groth_{\sigma\star_n \rho}(X).\]
Since the last descent of $\sigma$ is $l$, then
only the variables $x_1,\ldots,x_l$ appear in $\groth_{\sigma}(X)$. 
Hence because $l\leq t$, we have 
$r_{t}(\groth_{\sigma}(X))=\groth_{\sigma}(X)$ and so
\begin{equation}
\label{eqn:break}
  \groth_{\sigma}(X) \, r_{t}(\groth_{id\star_n \alpha}(X))
\!= \!r_t(\groth_{\sigma}(X)) \, r_{t}(\groth_{id\star_n \alpha}(X))
\!=\! r_{t}\big(\groth_{\sigma}(X)\, \groth_{id\star_n \alpha}(X)\big)
\!=\! r_{t}(\groth_{\sigma\star_n \alpha}(X)). 
\end{equation}
We conclude by applying the first conclusion to 
$\gamma=\sigma\star_n \alpha$ and observing that $\ell(\sigma\star_n \alpha)=
\ell(\sigma)+\ell(\alpha)$.
\end{proof}

\noindent 
\emph{Proof of Corollary~1:}
Since $r_t(\groth_{id\star_n\rho}(X))=\groth_{\rho}(X)$
~\cite{fomin.kirillov:1}, it follows from our above discussions that 
the hypotheses of Theorem~\ref{thm:main} hold.\qed

\noindent
{\bf Extensions to equivariant cohomology.}
The structure constants in the equivariant cohomology ring $H^*_T(\Fln)$ 
are polynomials in a second set
of variables $\{y_1,y_2,\ldots,\}$, and are known to have a positive 
expansion in $\{y_{i+1} - y_i\}$ (as proven in \cite{graham}).
It does not seem easy to extend our techniques to apply to this richer problem.
While the transition formula does have an equivariant
extension, and one can state an equivariant truncation formula,
this formula involves the $\{y_i\}$ individually rather than as differences.
In trying to group them into differences one leaves the realm
of subtraction-free formulae.

\end{document}